\documentclass[12pt]{article}
\usepackage{amsmath}
\usepackage{amsfonts}

\advance \textheight by \topskip
\advance \textheight by 1pt
\makeatother
\def\bea{\begin{eqnarray}}
\def\ena{\end{eqnarray}}

\def\lar{\longrightarrow}
\def\non{\nonumber}

\def\gr{\mathrm{gr}}
\def\dim{\mathrm{dim}}
\def\ch{\mathrm{ch}\,}
\def\ker{\mathrm{Ker}}

\def\tilp{\tilde{p}}

\def\tilp0{\tilde{p}_0}

\def\lar{\longrightarrow}

\def\hpsi{\widehat{\psi}}
\def\tpsi{\tilde{\psi}}

\newcommand{\bc}[2]{
\left(
\begin{array}{c}{#1}\\{#2}\end{array}
\right)}

\newcommand{\qed}{\hbox{\rule[-2pt]{3pt}{6pt}}}
\newtheorem{prop}{Proposition}
\newtheorem{theorem}{Theorem}

\newtheorem{lemma}{Lemma}

\title{
Baker-Akhiezer Modules on the Intersections of Shifted Theta Divisors
}

\author{
Koji Cho\thanks{
e-mail: cho@math.kyushu-u.ac.jp}\\
Department of Mathematics,\\
 Kyushu University\\
\\
Andrey Mironov\thanks{
e-mail: mironov@math.nsc.ru}\\
Sobolev Institute of Mathematics\\
and Novosibirsk State University\\
 \\
Atsushi Nakayashiki\thanks{
e-mail: 6vertex@math.kyushu-u.ac.jp}\\
Department of Mathematics,\\
 Kyushu University\\
}

\date{}
\begin{document}
\maketitle

\vskip2mm
\centerline{}
\vskip15mm

\begin{abstract}
\noindent
The restriction, on the spectral variables, of the Baker-Akhiezer (BA) 
module of a $g$-dimensional principally polarized abelian variety 
with the non-singular theta divisor 
to an intersection of shifted theta divisors is studied.
It is shown that the restriction to a $k$-dimensional variety
becomes a free module over the ring of differential operators
in $k$ variables. The remaining $g-k$ derivations define evolution
equations for generators of the BA-module.
As a corollary new examples of commutative 
ring of partial differential operators with matrix coefficients and
their non-trivial evolution equations are obtained.

\end{abstract}
\clearpage

\clearpage

\section{Introduction}
The Baker-Akhiezer(BA) module was introduced in \cite{N1,N2} in order to extend the theory of the BA function due to Krichever \cite{Kr} to higher dimensions.
It is a geometric counterpart of the ${\cal D}$-module generated by the 
wave operator in Sato's theory of KP-hierarchy and universal Grassmann manifold.

A fundamental example of the BA function is a function on an elliptic 
curve of the form
\bea
&&
\varphi(z;x)=\frac{\sigma(z+x)}{\sigma(z)\sigma(x)}e^{-x\zeta(z)},
\non
\ena
where $\sigma(z)$, $\zeta(z)$ are Weierstrass' sigma and zeta functions.
The corresponding BA-module is the ${\cal D}$-module generated by
$\varphi(z;x)$:
\bea
&&
M={\cal D}\varphi(z;x)=\sum_{n=0}^\infty {\cal O}\partial_x^n\varphi(z;x),
\non
\ena
where ${\cal O}$ is a suitable ring of functions such as the convergent 
power series ring, its quotient field etc. and ${\cal D}={\cal O}[\partial_x]$
is the ring of differential operators in $x$ with the coefficients in 
${\cal O}$. It is a rank one free module over ${\cal D}$.

Let $A={\mathbb C}[\wp(z),\wp'(z)]$ be the affine ring of the elliptic curve.
An important property of the BA-module is that it is not only a ${\cal D}$ 
module but also an $A$-module.
As a consequence $A$ is embedded in ${\cal D}$ as a commutative subring.

Similarly, in the case of genus $g$ algebraic curves,
 the BA-module becomes a $({\cal D}_g,A)$-bimodule, where  
${\cal D}_g={\cal O}[\partial_1,...,\partial_g]$ is the ring
of differential operators in $g$ variables and $A$ is 
the affine ring of the curve. It becomes a rank one
free module over the subring ${\cal D}_1={\cal O}[\partial_1]$ 
of operators in one variable and the affine ring $A$ is embedded
in ${\cal D}_1$. The action of the commuting derivations 
$\partial_2,...,\partial_g$ specifies evolution equations of the BA-module,
or the deformation of the image of $A$ in ${\cal D}$. In this way
solutions of integrable nonlinear equations such as KP equation , KdV equation
are constructed \cite{Kr}.

In \cite{N1} the BA-module of a $g$ dimensional principally polarized
Abelian variety $(X,\Theta)$ with a non-singular $\Theta$ is studied. 
It is proved that BA-module becomes a free
${\cal D}$-module of rank $g!$, where ${\cal D}$ is the ring of differential
operators in $g$ variables. Consequently the affine ring $A$ of 
$X\backslash\Theta$ is embedded in the ring $Mat(g!,{\cal D})$ of
differential operators with the coefficients in $g!\times g!$ matrices.
However in this case there is no non-trivial deformation. To have 
deformations it is necessary to consider the BA-module of
polarized subvarieties of $(X,\Theta)$. 

We consider an intersection $Y^k$ of shifted theta divisors as a subvariety
of $X$ and the intersection $Q^k$ of it with the theta divisor as a divisor,
where $k$ denotes the codimension of $Y^k$ in $X$.
We show that the restriction of the BA-module of $(X,\Theta)$ to $Y^k$
 is a free ${\cal D}_{g-k}$ module of rank $g!$, where ${\cal D}_i$ is
the ring of differential operators in $i$ variables. As a by-product
we have the embedding of the affine ring of $Y^k\backslash Q^k$ in
$Mat(g!,{\cal D}_{g-k})$ and $k$ commuting derivations 
which specify the deformation of the image of it. 

The simplest case of $g=3$ and
$k=1$ is studied in \cite{N2}. The case of intersections of more
general divisors are studied in \cite{Mir}. Unfortunately the proof of
the freeness is incomplete in that paper. 
Other examples of BA-modules which have non-trivial deformations 
are studied in \cite{R}.

The plan of the paper is as follows. In section 2 the definition of
BA-modules and the main 
result are given. The combinatorial properties of the restriction
of the BA-module is studied in section 3. 
It is shown that the character of the associated graded 
module of the BA-module coincides with that
of the free module. It means that as far as the dimension
is concerned the restriction of the BA-module becomes a free module.
In section 4 the proof of the main theorem is given based on the result
of section 3.

\section{Results}
Let $\Omega$ be a point of the Siegel upper half space of degree $g$,
$X$ the corresponding principally polarized Abelian variety
\bea
&&
X={\mathbb C}^g/\Gamma,
\quad
\Gamma={\mathbb Z}^g+\Omega{\mathbb Z}^g,
\non
\ena
$\theta_{a,b}(z)$ Riemann's theta function with the charcteristic ${}^t(a,b)$,$a,b\in {\mathbf R}^g$
\bea
&&
\theta_{a,b}(z)=\theta_{a,b}(z,\Omega)=\sum_{n\in {\mathbb Z}^g}\exp(\pi i{}^t (n+a)\Omega (n+a)
+2\pi i{}^t (n+a)(z+b)),
\non
\ena
and $\Theta$ the theta divisor on $X$ defined by the zero set of $\theta(z)=\theta_{0,0}(z)$.
For $c\in {\mathbb C}^g$, ${\cal L}_c$ denotes the holomorphic flat 
line bundle on $X$ which has $\theta(z+c)/\theta(z)$ as a meromorphic section.
A meromorphic section of ${\cal L}_c$ can be considered as a meromorphic 
function $f(z)$ on ${\mathbb C}^g$ which satisfies
\bea
&&
f(z+m+\Omega n)=\exp(-2\pi i{}^t nc)f(z),
\quad
m,n\in {\mathbb Z}^g.
\label{Lc}
\ena
We denote $L_c(n)$ the space of meromorphic sections of ${\cal L}_c$
whose poles are only on $\Theta$ of order at most $n$. 
It is known that $\dim\,L_c(n)=n^g$ and a linear basis is given
by the functions $f_{n,a}(z+\frac{c}{n})/\theta(z)^n$, 
$a\in {\mathbb Z}^g/n{\mathbb Z}^g$,\cite{Mum2}
\bea
&&
f_{n,a}(z)=\theta_{\frac{a}{n},0}(nz,n\Omega).
\non
\ena

We define
\bea
&&
L_c=\cup_{n=0}^\infty L_c(n).
\non
\ena
The subspaces $L_c(n)$ define an increasing filtration of $L_c$.

Set
\bea
&&
\zeta_i(z)=\frac{\partial}{\partial z_i}\log\,\theta(z).
\non
\ena
It satisfies, for $m,n\in {\mathbb Z}^g$,
\bea
&&
\zeta_j(z+m+\Omega n)=\zeta_j(z)-2\pi i n_j.
\non
\ena

Let ${\cal O}$ be the convergent power series ring in $x=(x_1,...,x_g)$,
${\cal K}$ its quotient field, 
$\partial_{x_i}=\partial/\partial x_i$,
${\cal D}={\cal K}[\partial_{x_1},...,\partial_{x_g}]$ the ring of differential
operators with the coefficients in ${\cal K}$, ${\cal D}(n)=
\{\sum_{|\alpha|\leq n}a_\alpha \partial_x^\alpha\in {\cal D}\}$ the differential operators of order at most $n$ and $\gr\,{\cal D}=
\oplus {\cal D}(n)/{\cal D}(n-1)$ the commutative ring of principal symbols.

In general, for a module with an increasing filtration $M=\cup_{n} M(n)$, the
associated graded module is defined by
\bea
&&
\gr\,M=\oplus_{n}\,\gr_n\,M,
\qquad
\gr_n\,M=M(n)/M(n-1).
\non
\ena

Let
\bea
&&
M_c(n)=\sum_{a\in {\mathbb Z}^g/n{\mathbb Z}^g}
 {\cal K}\frac{f_{n,a}(z+\frac{c+x}{n})}{\theta(z)^n}
e^{-\sum_{i=1}^g x_i\zeta_i(z)},
\non
\\
&&
M_c=\cup_{n=0}^\infty M_c(n).
\non
\ena
The space $M_c(n)$ is a $n^g$ dimensional vector space over ${\cal K}$ and 
the subspaces $\{M_c(n)\}$ specify an increasing filtration of $M_c$.
Then, for $c\notin \Gamma$,
\bea
&&
\dim_{{\cal K}}\gr_n\,M_c=\dim_{{\cal K}}\gr_n\,L_c=n^g-(n-1)^g,
\quad
n\geq 1.
\label{dim-grn}
\ena

As a function of the variables $z$ any element of $M_c$ satisfies the equation
(\ref{Lc}). The differentiation in $x_i$ preserves this equation.
 Thus $M_c$ becomes a ${\cal D}$-module.
 It is introduced in
\cite{N1} and called the Baker-Akhiezer(BA) module of $(X,\Theta)$. Since
\bea
&&
\partial_{x_i} M_c(n)\subset M_c(n+1),
\non
\ena
$\gr\, M_c$ becomes a $\gr\,{\cal D}$-module.

Let $A$ be the affine ring of $X\backslash\Theta$. Analytically it is
described as
\bea
&&
A=\{\frac{F(z)}{\theta(z)^n}\,\vert\, \text{$F(z)$ is holomorphic on ${\mathbb C}^g$,}\quad \frac{F}{\theta^n}(z+\gamma)=\frac{F}{\theta^n}(z)\quad
\text{for any $\gamma\in \Gamma$}\,\}.
\non
\ena
If $f(z)$ satisfies (\ref{Lc}) and $a(z)\in A$ then $f(z)a(z)$ satisfies
(\ref{Lc}). Therefore $L_c$ is an $A$-module. Consequently $M_c$ becomes
an $A$-module whose action obviously commutes with that of ${\cal D}$.
\vskip2mm

\noindent
{\bf Remark} In \cite{N1} the BA module is defined globally
on the dual abelian variety $\widehat{X}$ of $X$ as a sheaf of 
${\cal D}_{\widehat{X}}$-module by Fourier-Mukai transform, 
 $M_c$ being a scalar extension of the stalk at the point specified
by $c$.
\vskip2mm

Let $r_i$, $i\geq 1$, be integers defined by
\bea
&&
r_i=i^g-(i-1)^g-\sum_{j=1}^{i-1}r_j\bc{g+i-j-1}{g-1}.
\non
\ena
They satisfy $r_1=1$, $r_i=r_{g-i}$, $r_i=0$ for $i>g$ and $\sum_{i=1}^g r_i=g!$\cite{N1}.

The following theorem is proved in \cite{N1}.

\begin{theorem}\label{th1}\cite{N1}
Suppose that $\Theta$ is non-singular and 
$c\notin {\mathbb Z}^g+\Omega {\mathbb Z}^g$.
Then $\gr\,M_c$ is a free $\gr\,{\cal D}$-module of rank $g!$ and $M_c$ is a free ${\cal D}$-module of rank $g!$. More precisely
\bea
\gr\,M_c&=&\oplus_{i=1}^g\oplus_{j=1}^{r_i} (\gr\,{\cal D}) \psi_{ij},
\quad
\psi_{ij}\in \gr_i\,M_c,
\non
\\
M_c&=&\oplus_{i=1}^g\oplus_{j=1}^{r_i} {\cal D} \phi_{ij},
\quad
\phi_{ij}\in M_c(i),
\ena
where $\psi_{ij}$ is the projection
of $\phi_{ij}$ in $\gr_i\,M_c$.
\end{theorem}

Let $\Phi=(\phi_{ij})$ be the column vector of dimension $g!$
 and $\text{Mat}(m,{\cal D})$ the ring of $m\times m$ matrices
with the components in ${\cal D}$.  Since $M_c$ is an $A$-module,
 for $a\in A$, there exists
an element $\ell(a)\in \text{Mat}(g!,{\cal D})$ such that
\bea
&&
a\Phi=\ell(a)\Phi.
\non
\ena
It defines an embedding of $A$ into $\text{Mat}(g!,{\cal D})$ as a ring.
Thus $A$ is realized as a commutative subring of the ring 
differential operators in $g$ variables with matrix coefficients.

In this paper we shall extend Theorem \ref{th1} to the BA-modules on the
intersections of shifted theta divisors.

For $a\in {\mathbb C}^g$ we set
\bea
&&
\Theta_a=\{\theta(z-a)=0\}\subset X.
\non
\ena
Take $a_1,...,a_{g-1}\in {\mathbb C}^g$. 
Beginning from $(Y^0,Q^0)=(X,\Theta)$ we define $(Y^k,Q^k)$ for $k\geq 1$
by
\bea
&&
Y^{k}=\Theta_{a_1}\cap\cdots\cap \Theta_{a_k},
\non
\\
&&
Q^{k}=Y^k\cap \Theta.
\non
\ena
We assume that, for any $k\leq g$, 
$\Theta_{a_1}$, ..., $\Theta_{a_k}$ intersect transversally and so does
$\Theta$, $\Theta_{a_1}$, ..., $\Theta_{a_k}$.
It means in particular that $Y^k$ and $Q^k$ are non-singular subvarieties of $X$ of dimensions $g-k$ and $g-k-1$ respectively. It can be shown that for generic $a_1$,...,$a_{g-1}$ the assumption is satisfied.

We denote the restriction of ${\cal L}_c$ to $Y^k$ by the same symbol
for simplicity. 
Let ${\cal L}_{c}(nQ^k)$ be the sheaf of germs of meromorphic sections
of ${\cal L}_c$ 
on $Y^k$ with poles only on $Q^k$ of order at most $n$. 
 We set 
\bea
&&
L^k_c(n)=H^0(Y^k,{\cal L}_c(nQ^k)).
\non
\ena
We define the space $\tilde{M}^k_c(n)$ as the set of functions of the form
\bea
&&
f(z;x)e^{-\sum_{i=1}^g x_i \zeta_i(z)\vert_{Y^k}},
\non
\ena
where $f(z;x)$ satisfies the following conditions. 

There is an open neighborhood $U$
 of $0\in{\mathbb C}^g$, which can depend on $f$, with the following properties.
\vskip2mm
\noindent
(i) For each $x\in U$ $f(z;x)$ belongs to $L^k_{c+x}(n)$ as a function of $z$.
\vskip2mm
\noindent
(ii) As a function of $x$ $f(z;x)$ is analytic on $U$.
\vskip2mm

It is obvious that ${\tilde M}_c^k(n)$ is an ${\cal O}$-module.

\begin{lemma}\label{lem-1} Let $k$ satisfy $0\leq k\leq g-1$. 
Suppose that
$c+\sum _{i\in I} a_i\notin \Gamma$ for any subset $I$ of $\{1,...,k\}$.
Then
\vskip2mm
\noindent
(i) $H^i(Y^k,{\cal L}_c(nQ^k))=0$, $i\neq 0,g-k$, $n\in {\mathbb Z}$.
\vskip2mm
\noindent
(ii) The restriction map 
$
L^{k-1}_c(n) \lar L^{k}_{c}(n)
$
is surjective for any $n\in {\mathbb Z}$.
\end{lemma}
\vskip2mm
\noindent
{\it Proof.} 
We have the following exact sequence and isomorphism:
\bea
&&
0\lar {\cal L}_c(nQ^{k-1}-Y^{k}) \lar {\cal L}_c(nQ^{k-1}) \lar 
{\cal L}_c(nQ^{k}) \lar 0,
\label{exact-1}
\\
&&
{\cal L}_c(nQ^{k-1}-Y^{k})\simeq {\cal L}_{c+a_{k}}((n-1)Q^{k-1}).
\label{isom}
\ena
The isomorphism (\ref{isom}) follows from
${\cal O}_X(-\Theta_a)\simeq {\cal O}_X(-\Theta)\otimes {\cal L}_a$.
The assertion (i) can be proved by induction on $k$ using the cohomology 
sequence of $(\ref{exact-1})$
and the vanishing \cite{Mum1}:
\bea
&&
H^i(X,{\cal L}_c(n\Theta))=0,
\non
\ena
for $i\geq 1$, $n\geq 1$ or $i\geq 0$, $n=0$ or $i\neq g$, $n<0$.
The assertion (ii) follows from the cohomology sequence of (\ref{exact-1})
and (i).
\qed
\vskip2mm

\begin{lemma}\label{lem-1-2}
Assume the same conditions as in Lemma \ref{lem-1}. 
Then $\dim\,L^k_c(n)$ does not
depend on $c$ and satisfies
\bea
&&
\dim\,L^{k}_c(n)=\dim\,L^{k-1}_c(n)-\dim\,L^{k-1}_c(n-1).
\label{recursion}
\ena
\end{lemma}
\vskip2mm
\noindent
{\it Proof.} The lemma can easily proved by induction using 
the cohomology sequence 
of (\ref{exact-1}),
the isomorphism (\ref{isom}), (i) of Lemma \ref{lem-1} 
and the fact that $\dim\,L^0_c(n)=n^g$ for $n\geq 0$ if $c\notin \Gamma$.
\qed
\vskip2mm

\noindent
{\bf Example.} $\dim\,L^1_c(n)=n^g-(n-1)^g$ for $n\geq 1$.
\vskip2mm
\noindent
$\dim\,L^2_c(n)=n^g-2(n-1)^g+(n-2)^g$ for $n\geq 2$ and $\dim\,L^2_c(1)=1$.

\begin{lemma}\label{lem-2}
Assume the same conditions as in Lemma \ref{lem-1-2}.
Then ${\tilde M}^k_c$ is a free ${\cal O}$-module of rank 
$\dim\,L^k_c(n)$.
\end{lemma}
\vskip2mm
\noindent
{\it Proof.} 
We take some basis $\{f_i(z)\}$ of
$L^k_c(n)$ and lift it to $\{f_i(z;x)\}$ such that
the conditions (i),(ii) are satisfied and $f_i(z;0)=f_i(z)$.
 Then it gives a ${\cal O}$-free basis of 
${\tilde M}_c^k(n)$. Such analytic lift can be constructed among
the functions $\{f_{n,a}(z+\frac{c+x}{n})\}$ restricted to $Y^k$
since the restriction map $L^0_{c+x}(n)\rightarrow L^k_{c+x}(n)$ 
is surjective for $x$ sufficiently close to $0\in {\mathbb C}^g$ by (ii)
of Lemma \ref{lem-1}.
\qed

We set
\bea
&&
M^k_c(n)={\cal K}\otimes_{\cal O}{\tilde M}_c^k,
\quad
M^k_c=\cup M^k_c(n).
\non
\ena
The set of subspaces $\{M_c^k(n)\}$ 
defines an increasing filtration on $M_c^k$.

Let $\pi_k:M_c^k\lar M_c^{k+1}$ be the restriction map. It is surjective for
$k\geq 0$ as shown above. 
In particular $\pi_{0k}:=\pi_k\pi_{k-1}\cdots \pi_0: M_c\lar M_c^k$ 
is surjective. 
Thus $M^k_c$ can be directly described as the restriction of 
$M_c$ to $Y^k$ with respect to the $z$ variables:
\bea
&&
M^k_c=M_c\vert_{Y^k}.
\non
\ena
It is obvious that the restriction in $z$ variables commutes with the action
of $\partial_{x_i}$. Therefore $M^k_c$ becomes a ${\cal D}$-module.
Moreover the action of $\partial_{x_i}$ satisfies 
$\partial_{x_i}M^k_c(n)\subset M^k_c(n+1)$. Thus $\gr\,M^k_c$ becomes a
$\gr\,{\cal D}$ module. 
The main result of this paper is

\begin{theorem}\label{th2}
Suppose that $\Theta$ is non-singular and 
$c+\sum_{i\in I} a_i\notin \Gamma$ for any subset $I$ of
$\{1,2,...,g-1\}$. Then
there exists a set of linear independent vector fields 
$D_i=\sum_{j=1}^g c_{ij}{\partial_{x_j}}$, $c_{ij}$ being constants, such that
the following properties are valid.
\vskip2mm
\noindent
(i) Let ${\cal D}_{i}={\cal K}[D_1,...,D_{i}]$. Then
$M^k_c$ is a free ${\cal D}_{g-k}$ module of rank $g!$. More precisely
it is of the form
\bea
&&
M^k_c=\oplus_{i=1}^g\oplus_{j=1}^{r_i}{\cal D}_{g-k}\phi^k_{ij},
\quad
\phi^k_{ij}\in M^k_c(i).
\non
\ena
\vskip2mm
\noindent
(ii) 
The module
$\gr\,M^k_c$ is a free $\gr\,{\cal D}_{g-k}$ module of rank $g!$. 
More precisely it is of the form
\bea
&&
\gr\,M^k_c=\oplus_{i=1}^g\oplus_{j=1}^{r_i}(\gr\,{\cal D}_{g-k})\psi^k_{ij},
\quad
\psi^k_{ij}\in\gr_i\,M^k_c,
\non
\ena
where $\psi^k_{ij}$ is the projection of $\phi^k_{ij}$
in $\gr_i\,M^k_c$ and the filtration of ${\cal D}_{g-k}$ is specified by
${\cal D}_k(n)={\cal D}_k\cap {\cal D}(n)$.
\end{theorem}

Let $\Phi^k=(\phi^k_{ij})$ be the column vector
and $A^{g-k}$ be the affine ring of $Y^k\backslash Q^k$.
We have $\Phi^0=\Phi$, $A^0=A$. The ring $A^{g-k}$ acts on
$M^k_c$. Thus for any $a\in A^{g-k}$
there exists a unique operator $\ell^k(a)\in \text{Mat}(g!,{\cal D}_{g-k})$
such that
\bea
&&
a\Phi^k=\ell^k(a)\Phi^k.
\non
\ena
It defines an embedding of $A^{g-k}$ in $\text{Mat}(g!,{\cal D}_{g-k})$.
Since $M^k_c$ is a ${\cal D}$-module and ${\cal D}={\cal D}_g$,
for $D_i$ with $g-k+1\leq i\leq g$, there exists a unique operator
$B^k_i\in \text{Mat}(g!,{\cal D}_{g-k})$ such that
\bea
&&
D_i\Phi^k=B^k_i\Phi^k.
\non
\ena
Those operators satisfy, for any $a,b,i,j$,
\bea
&&
[\ell^k(a),\ell^k(b)]=0, 
\non
\\
&&
[D_i-B^k_i,D_{j}-B^k_{j}]=0,
\non
\\
&&
[D_i-B^k_i,\ell^k(a)]=0.
\non
\ena

\section{Combinatorial freeness}
For a graded ${\cal K}$-vector space
 $V=\oplus_{n\in {\mathbf Z}} V_n$ such that
each $V_n$ is finite dimensional we define the character
$\ch V$ by
\bea
&&
\ch V=\sum (\text{dim}_{{\cal K}} V_n)t^n.
\non
\ena
 Obviously
\bea
&&
\ch \gr\,{\cal D}_j=\frac{1}{(1-t)^j}.
\non
\ena
We have
\bea
\ch \gr\,M^0_c&=&\sum_{n=1}^\infty \left(n^g-(n-1)^g\right)t^n
\non
\\
&=&(1-t)\left(t\frac{d}{dt}\right)^g(1-t)^{-1}.
\non
\ena

Let $\{\psi_{ij}\}$ be a $\gr\,{\cal D}$-free basis of $\gr\,M_c$ 
as in Theorem 1 and $F=\oplus {\cal K}\psi_{ij}$ the subspace
of $\gr\,M_c$. By Theorem \ref{th1} 
$\gr\,M_c\simeq (\gr\,{\cal D}_g)\otimes F$.
The module $F$ naturally inherits a grading from $\gr\,M_c$. Then
\bea
&&
\ch\gr\,M_c=(\ch\,\gr\,{\cal D}_g)\cdot \gr\,F=(1-t)^{-g}\sum_{i=1}^g r_i t^i.
\non
\ena

\vskip2mm
\noindent
{\bf Example} For $g=1,2,3,4$ $\ch \gr\,M_c$ is given by
\bea
&&
\frac{t}{(1-t)},
\quad
\frac{t+t^2}{(1-t)^2},
\quad
\frac{t+4t^2+t^3}{(1-t)^3},
\quad
\frac{t+11t^2+11t^3+t^4}{(1-t)^4}.
\non
\ena
\vskip2mm

\begin{lemma}\label{lem-3-0} Assume the same conditions for $c,a_1,...a_{g-1}$ 
as in Theorem \ref{th2}. Then
\vskip2mm
\noindent
(i) $\dim_{\cal K}\, \gr_n\,M_c^{j+1}=\dim_{\cal K}\,\gr_n\,M_c^{j}-
\dim_{\cal K}\,\gr_{n-1}\,M_c^{j}$.
\vskip2mm
\noindent
(ii)
$\displaystyle{
\ch\gr\,\,M^{j+1}_c=(1-t)\ch\gr\,\,M^j_c}$.
\end{lemma}
\vskip2mm
\noindent
{\it Proof.} The assertion (ii) follows from (i) and (i) follows from
Lemma \ref{lem-1-2} and \ref{lem-2}. 
\qed
\vskip2mm

By Lemma \ref{lem-3-0} we have
\bea
&&
\ch\gr\,\,M^j_c=(1-t)^{-(g-j)}\sum_{i=1}^g a_i t^i
=(\ch {\cal D}_{g-j})\cdot \ch F.
\label{char-identity}
\ena

\section{Proof of Theorem \ref{th2}}

Notice that (i) of Theorem \ref{th2} follows from (ii) of Theorem \ref{th2}.

We shall prove

\begin{prop}\label{prop1}
Assume the same conditions as in Theorem \ref{th2}.
Set $y^{(0)}=(y^{(0)}_1,...,y^{(0)}_g)=(x_1,...,x_g)$.
Then, for each $k\geq 1$ there exist  a linear change of the coordinates
from $y^{(k-1)}=(y^{(k-1)}_1,...,y^{(k-1)}_{g-k+1})$ to
$y^{(k)}=(y^{(k)}_1,...,y^{(k)}_{g-k+1})$ and 
$\psi_{ij}^{k}\in \gr_i\,M_c^k$, $1\leq i\leq g$, $1\leq j\leq r_i$ such that
the following properties hold. 
Let 
${\cal D}_{g-k}={\cal K}[\partial_{y^{(k)}_1},...,\partial_{y^{(k)}_{g-k}}]$, 
${\cal D}_{g-k}(n)={\cal D}_{g-k}\cap {\cal D}(n)$ and
$\xi^{(k)}_i$ the image of $\partial_{y^{(k)}_i}$ in $\gr_1\,{\cal D}$.
Then
\bea
&&
\xi^{(k)}_{g-k+1}\psi_{ij}^k\in \sum_{i'j'}(\gr\,{\cal D}_{g-k})\psi^k_{i'j'},
\non
\\
&&
\gr\,M_c^k=\oplus_{i=1}^g\oplus_{j=1}^{r_i}(\gr\,{\cal D}_{g-k})\psi^k_{ij}.
\non
\ena
\end{prop}

If we define, for $1\leq k\leq g$, $D_{g-k+1}=\partial_{y^{(k)}_{g-k+1}}$, then
Theorem \ref{th2} (ii) follows from this proposition.

We prove the proposition by induction on $k$, 
 where the case of $k=0$ is established by Theorem \ref{th1}. 
We assume that the proposition is valid for $k$ if 
$c+\sum_{i\in I} a_i\notin \Gamma$ for any subset $I$ of $\{1,...,k\}$. 
We shall prove that
the proposition is true for $k+1$ if $c+\sum_{i\in I} a_i\notin \Gamma$ 
for any subset $I$ of $\{1,...,k+1\}$.

Let us set 
\bea
&&
\tpsi_{ij}^{k+1}=\psi_{ij}^{k}\vert_{Y^{k+1}}.
\non
\ena

\begin{lemma}\label{lem1}
For each $(ij)$ there exist a non-zero element $P_{ij}\in \gr_{N^{ij}}\,{\cal D}_{g-k}$ for some $N^{ij}\geq 0$ and  a linear change of the coordinates from 
$(y^{(k)}_1,...,y^{(k)}_{g-k})$ to $(y^{(k+1)}_1,...,y^{(k+1)}_{g-k})$
such that the following properties hold.
\vskip2mm
\noindent
(i) $
P_{ij}\tpsi_{ij}^{k+1}=0\quad \text{in $\gr\,M^{k+1}_c$}.
$
\vskip2mm
\noindent
(ii) Let $\xi_i=\xi^{(k+1)}_i$.  Then $P_{ij}$ is of the form
\bea
&&
P_{ij}=\xi_{g-k}^{N^{ij}}+\sum_{|\alpha|=N^{ij},\alpha_{g-k}< N^{ij}}
a_{ij;\alpha}\xi_1^{\alpha_1}\cdots\xi_{g-k}^{\alpha_{g-k}},
\non
\ena
where $\alpha=(\alpha_1,...,\alpha_{g-k})$, 
$|\alpha|=\sum_{i=1}^{g-k} \alpha_i$.
\end{lemma}
\vskip2mm
\noindent
{\it Proof.} (i) 
By (\ref{char-identity}) the dimension of $\gr_n\,M^{k+1}_c$ 
is a polynomial in $n$ of degree
$g-k-2$ for sufficiently large $n$. 
If there are no non-trivial linear relations among 
$\xi_1^{\alpha_1}\cdots\xi_{g-k}^{\alpha_{g-k}}
\tpsi^{k+1}_{ij}$, $\sum \alpha_l=n-i$
for any $n$, then $\dim_{\cal K} 
\left((\gr_{n-i}\,{\cal D}_{g-k})\tpsi_{ij}^{k+1}\right)=
\dim_{\cal K}\gr_{n-i}{\cal D}_{g-k}$ and it is a polynomial in $n$ of degree
$g-k-1$ for sufficiently large $n$. Thus there should be a relation as in the assertion.

\vskip2mm
\noindent
(ii) Let us write 
\bea
&&
P_{ij}=\sum_{\sum \alpha_l=N^{ij}} q_{\alpha_1,...,\alpha_{g-k}} 
\xi_1^{\alpha_1}\cdots \xi_{g-k}^{\alpha_{g-k}},
\non
\ena
Changing the name of the variables if necessary one can assume that
$q_{\alpha_1,...,\alpha_{g-k}}\neq 0$ for some $(\alpha_1,...,\alpha_{g-k})$
with $\alpha_{g-k}\neq 0$. If $q_{0,...,0,N^{ij}}\neq 0$ then we get the
desired element by dividing $P_{ij}$ by $q_{0,...,0,N^{ij}}\neq 0$.
If this is not the case, we consider the change of the variables of the
form
\bea
&&
\xi^{(k)}_i=\sum_{l=i}^{g-k}c_{i,l}\xi^{(k+1)}_l.
\non
\ena
Let $c_i=c_{i,g-k}$.
Then in the resulting expression of $P_{ij}$ the coefficient of
$(\xi^{(k+1)}_{g-k})^{N^{ij}}$ is
\bea
&&
\sum_{|\alpha|=N^{ij}} c_{1}^{\alpha_1}\cdots c_{g-k}^{\alpha_{g-k}}q_{\alpha_1,...,\alpha_{g-k}}.
\label{coeff}
\ena
This is a non-zero homogeneous polynomial in $c_1,...,c_{g-k}$. 
Thus it is non-zero on a non-empty open subset of ${\mathbf C}^{g-k}$.
Take a point of it, make a change of the coordinates and dividing $P_{ij}$
by (\ref{coeff}) we get a desired result.
\qed

Let 
\bea
&&
\gr\,\pi_k:\gr\,M^k_c\lar \gr\,M^{k+1}_c
\non
\ena
be the restriction map induced by $\pi_k$ and $K^k=\oplus K^k_n$
 the kernel of $\gr\,\pi_k$. 
Since $\gr\,\pi_k$ is a homomorphism of $\gr\,{\cal D}$-modules,
 $K^k$ is a $\gr\, {\cal D}$ submodule of $\gr\,M^k_c$. 
We denote by ${\tilde K}^k=\oplus {\tilde K}^k_n$ the $\gr\, {\cal D}$ module
obtained from $K^k$ by shifting the grading by $-1$, that is, 
${\tilde K}^k_n=K^k_{n+1}$.

\begin{lemma}\label{lem-5} The map
\bea
\gr_n\, M^k_{c+a_{k+1}}&\lar& K^k_{n+1}
\non
\\
\phi(z)&\mapsto& \frac{\theta(z-a_{k+1})}{\theta(z)}\phi(z),
\label{map}
\ena
gives an isomorphism of $\gr\,M^k_{c+a_{k+1}}$ and ${\tilde K}^k$ as 
 $\gr\,{\cal D}_{g-k}$-modules.
\end{lemma}
\vskip2mm
\noindent
{\it Proof.} We can assume $g-k\geq 2$. Using Lemma \ref{lem-1} (i), (\ref{isom}) and the cohomology sequence of (\ref{exact-1}) we have 
\bea
&&
\ker(\pi_k|_{M^k_c(n)})=
\frac{\theta(z-a_{k+1})}{\theta(z)}|_{Y^k}M^k_{c+a_{k+1}}(n-1).
\label{ker-pi-k}
\ena
Let $\gr_n \left({\cal L}_c(-m Y^{k+1})|_{Y^k}\right)$ be the sheaf
 on $Y^k$ defined by the exact sequence;
\bea
0\rightarrow {\cal L}_c((n-1)Q^k-mY^{k+1}) 
&\rightarrow&{\cal L}_c(nQ^k-mY^{k+1})
\non
\\
&\quad&
\quad
\rightarrow \gr_n \left({\cal L}_c(-mY^{k+1})|_{Y^k}\right) 
\rightarrow 0.
\label{gr-def}
\ena
Then one can easily verify that the following is an exact sequence,
\bea
&&
0\lar \gr_n\,\left({\cal L}_c(-Y^{k+1})|_{Y^k}\right)
 \lar \gr_n\, ({\cal L}_c|_{Y^k})
\lar \gr_n\, ({\cal L}_c|_{Y^{k+1}})\lar 0.
\label{gr-exact}
\ena
By the cohomology sequence of (\ref{gr-def}) with $m=0$ the natural map
\bea
&&
\gr_n\,L^k_c:=L^k_c(n)/L^k_c(n-1)\lar 
H^0\left(Y^k,\gr_n({\cal L}_c|_{Y^k})\right)
\label{gr-to-gr}
\ena
is always injective and becomes isomorphic if $g-k\geq 2$. Then we have
\bea
\ker(\gr_n\,L^k_c\rightarrow \gr_n\,L^{k+1}_c)
&\simeq&H^0\left(Y^k,\gr_n({\cal L}_c(-Y^{k+1})|_{Y^k})\right)
\non
\\
&\simeq&
\frac{H^0(Y^k,{\cal L}_c(nQ^k-Y^{k+1}))}{H^0(Y^k,{\cal L}_c((n-1)Q^k-Y^{k+1}))}
\non
\\
&\simeq&
\frac{\theta(z-a_{k+1})}{\theta(z)}|_{Y^k}
\gr_{n-1}\,L^k_{c+a_{k+1}},
\non
\ena

where we use the cohomology sequences of (\ref{gr-exact}), (\ref{gr-def}),
Lemma \ref{lem-1} (i), (\ref{isom}), (\ref{ker-pi-k}).
\qed
\vskip2mm

If $c+a_{k+1}+\sum_{i\in I} a_i\not\in \Gamma$ for any subset $I$ of
$\{1,..,k\}$ the induction hypothesis can be applied to
$M_{c+a_{k+1}}^k$.
By the lemma and the assumption of induction $K^k$ is described as
\bea
&&
K^k=\oplus_{i=1}^g\oplus_{j=1}^{r_i}(\gr\,{\cal D}_{g-k})\varphi_{ij},
\non
\\
&&
\varphi_{ij}=\frac{\theta(z-a_{k+1})}{\theta(z)}
\psi^k_{ij}\vert_{c\rightarrow c+a_{k+1}}
\in K^k_{i+1}.
\non
\ena

Since $\varphi_{ij}\in \gr_{i+1}\,M_c^{k}$, it can be written as a linear
combination of $\{\psi^k_{ij}\}$ with the coefficients in $\gr\,{\cal D}_{g-k}$
as
\bea
&&
\varphi_{ij}=\sum Q_{ij:i'j'}\psi^k_{i'j'},
\quad
Q_{ij:i'j'}=\sum_{|\alpha|+i'=i+1}q^{\alpha}_{ij;i'j'}\xi^\alpha.
\label{eq-1}
\ena

By (i) of Lemma \ref{lem1} we have $P_{ij}\psi^k_{ij}\in K^k_{i+N^{ij}}$.
Thus it can be written as a linear combination of $\{\varphi_{ij}\}$ 
with the coefficients in $\gr\,{\cal D}_{g-k}$ as
\bea
&&
P_{ij}\psi^k_{ij}=\sum R_{ij;i'j'}\varphi_{i'j'},
\quad
R_{ij;i'j'}=\sum_{|\alpha|+i'+1=i+N^{ij}} r^\alpha_{ij;i'j'}\xi^\alpha.
\label{eq-2}
\ena
Composing these relations we get
\bea
&&
P_{ij}\psi^k_{ij}=\sum R_{ij;i'j'}Q_{i'j':i''j''}\psi^k_{i''j''}.
\label{eq-3}
\ena
In the matrix form it is written as
\bea
&&
P=RQ,
\label{eq-4}
\ena
where $P=(P_{ij})$ is the diagonal matrix and $R=(R_{ij;i'j'})$, 
$Q=(Q_{ij;i'j'})$ are $g!\times g!$ matrices.

We shall construct a basis $\{\hpsi^k_{ij}\}$ of $\gr\,M^k_c$ 
as a $\gr\,{\cal D}_{g-k}$-module modifying $\{\psi^k_{ij}\}$
appropriately such that they satisfy
\bea
&&
\xi_{g-k}\hpsi^k_{ij}\in \sum_{i'\leq i+1}{\cal K}[\xi_1,...,\xi_{g-k-1}]\hpsi^k_{lm}+K^k.
\label{aim1}
\ena
To this end we use the relation (\ref{eq-1}). Let us write 
it more explicitly as
\bea
&&
\varphi_{ij}=\sum q_{ij;i+1j'}\psi^k_{i+1,j'}
+ \sum_{l=1}^{g-k}\sum_{j'}q^l_{ij;ij'}\xi_l \psi^k_{ij'}
+\sum_{i'<i,|\alpha|+i'=i+1}q^\alpha_{ij;i'j'}\xi^\alpha\psi^k_{i'j'}.
\label{eq-5}
\ena
For the sake of simplicity we identify $q^l_{ij;ij'}$ with 
$q^{0,...,0,l}_{ij;ij'}$.

We construct $\{\hpsi^k_{ij}\}$ satisfying the property (\ref{aim1}) 
by induction on $i$. 

Let us consider the case $i=1$. Then the equation (\ref{eq-5}) becomes
\bea
&&
\varphi_{11}=\sum q_{11;2j}\psi^k_{2j}+\sum_{l=1}^{g-k}q_{11;11}^l\xi_l\psi^k_{11}.\non
\ena
We proceed by dividing the case.
\vskip2mm
\noindent
(i) the case $q^{g-k}_{11;11}\neq 0$. In this case (\ref{aim1}) holds for 
$(ij)=(11)$ by defining $\hpsi^k_{11}=\psi^k_{11}$ and 
$\hpsi^k_{2j}=\psi^k_{2j}$.
\vskip2mm
\noindent
(ii) the case $q^{g-k}_{11;11}=0$ and $q_{11;2j}\neq 0$ for some $j$.
In this case we modify $\psi^k_{2j}$ to $\hpsi^k_{2j}=\psi^k_{2j}-
\xi_{g-k}\psi^k_{11}$. Correspondingly the operators $P_{ij}$ is changed.
We can take the product $P_{2j}P_{11}$ as the operator for $\hpsi^k_{2j}$,
since $P_{2j}P_{11}\hpsi^k_{2j}=0$. 
We set $\hpsi^k_{2j'}=\psi^k_{2j'}$ for $j'\neq j$ and 
$\hpsi^k_{11}=\psi^k_{11}$.
Then (\ref{aim1}) holds for $(ij)=(11)$.
\vskip2mm
\noindent
(iii) the case $q^{g-k}_{11;11}=0$ and $q_{11;2j}= 0$ for all $j$.
This case is impossible. In fact, suppose that this is the case. We consider
the equation (\ref{eq-4}) modulo the ideal
\bea
&&
I=\sum_{l=1}^{g-k-1} (\gr\,{\cal D}_{g-k})\xi_l.
\non
\ena
Then the right hand side is degenerate while the left hand side 
is non-degenerate.
\vskip2mm

As a whole we have constructed $\hpsi^k_{ij}$, $i\leq2$ such that
(\ref{aim1}) holds for $(ij)=(11)$.

Assume that $\hpsi^k_{i'j'}$, $i'\leq i$, $1\leq j'\leq r_{i'}$, 
are constructed in such a way that (\ref{aim1}) holds for $(i'j')$, $i'<i$.
As a consequence of the change from $\{\psi^k_{i'j'}|i'\leq i\}$ to
$\{\hpsi^k_{i'j'}|i'\leq i\}$ the matrices $P$, $Q$, $R$ may be changed.
However their properties that the equation (\ref{eq-4}) holds and 
$P$ is non-degenerate remain valid. Therefore we use the same symbol
$P$, $Q$, $R$ and their components in the argument below for the sake
of simplicity.

Using the relation (\ref{aim1}) for $\hpsi^k_{i'j'}$ for $i'<i$, the equation
(\ref{eq-1}) for $\varphi_{ij}$ can be written as
\bea
&&
{\tilde \varphi}_{ij}=\sum q_{ij;i+1j'}\psi^k_{i+1,j'}
+\sum_{l=1}^{g-k}\sum_{j'=1}^{r_i}{\tilde q}^l_{ij;ij'}\xi_l\hpsi^k_{ij'}
+\sum_{i'<i,|\alpha|+i'=i+1,\alpha_{g-k}=0} {\tilde q}^{\alpha}_{ij;i'j'}
\xi^\alpha\hpsi^k_{i'j'},
\non
\ena
for some ${\tilde \varphi}_{ij}\in K^k_{i+1}$ and some 
${\tilde q}^{\alpha}_{ij;i'j'}$, $i'\leq i$. Here 
$q^{\alpha}_{ij;i'j'}$ changes to ${\tilde q}^{\alpha}_{ij;i'j'}$ as a consequence of the use of (\ref{aim1}). 
Notice that to use (\ref{aim1}), as to the effect on the matrix $Q$, is 
 to make fundamental transformations in rows of $Q$.

We have
\bea
&&
{\tilde \varphi}_{ij}=\sum q_{ij;i+1j'}\psi^k_{i+1,j'}
+\sum_{j'=1}^{r_i}{\tilde q}_{ij;ij'}^{g-k}\xi_{g-k}\psi^k_{ij'}
\quad\text{mod. $I\gr\,M^k_c$.}
\non
\ena
The rank of the $r_i\times (r_i+r_{i+1})$ matrix
\bea
&&
\left((q_{ij;i+1j'})_{1\leq j'\leq r_{i+1}},
({\tilde q}^{g-k}_{ij;ij'})_{1\leq j'\leq r_{i}}\right)_{1\leq j\leq r_i}
\non
\ena
is maximal. For, otherwise it contradicts the non-degeneracy of the matrix $P$.
Thus, as in the case of $i=1$, modifying $\psi^k_{i+1j'}$ to 
$\hpsi^k_{i+1j'}=\psi^k_{i+1j'}-\xi_{g-k}\psi^k_{ij''}$ for some $j'$ and $j''$
if necessary, we get $\{\hpsi^k_{i',j'}\,|i'\leq i+1\}$ such that (\ref{aim1})
holds for $(i',j')$, $i'\leq i$.
Notice that, in the last step $i=g$, $\tilde{q}_{g+11;g1}^{g-k}\neq 0$ 
is automatic.  

Thus a $\gr\,{\cal D}_{g-k}$-free basis $\{\hpsi^k_{ij}\}$ of $\gr\,M^k_c$
satisfying the condition (\ref{aim1}) is constructed. Set
\bea
&&
\psi^{k+1}_{ij}=\hpsi^k_{ij}\vert_{Y^{k+1}}.
\non
\ena
Then (\ref{aim1}) implies that
\bea
&&
\xi_{g-k}\psi^{k+1}_{ij}\in \sum \gr\,{\cal D}_{g-k-1}\psi^{k+1}_{i'j'}.
\label{reduction}
\ena

\begin{lemma} 
If $c+\sum_{i\in I} a_i\notin \Gamma$ for any subset $I$ of $\{1,...,k\}$, 
the restriction map $\gr\,\pi_k$ is surjective.
\end{lemma}
\vskip2mm
\noindent
{\it Proof.}
The lemma can be proved by a similar argument to the proof of 
Lemma \ref{lem-5}.
\qed
\vskip2mm

By the lemma and the assumption of induction on $\gr\,M^{k}_c$, $\gr\,M^{k+1}_c$ is generated by $\{\psi^{k+1}_{ij}\}$
over $\gr\,{\cal D}_{g-k}$. Then (\ref{reduction}) implies that
\bea
&&
\gr\,M^{k+1}_c=\sum (\gr\,{\cal D}_{g-k-1})\psi^{k+1}_{ij}.
\non
\ena
It follows from (\ref{char-identity}) that the sum of the right hand side is a 
direct sum. This completes the proof of Proposition \ref{prop1}.
\qed

\vskip5mm
\noindent
{\large {\bf Acknowledgement}} 
\vskip3mm
\noindent
This research is supported by Grant-in-Aid for Scientific Research (B) 
17340048.

\end{document}